\documentclass[a4,11pt]{amsart}
\usepackage{amssymb,amsmath,amsthm,txfonts}
\usepackage{cite}

\usepackage{color}

\newtheorem{thm}{Theorem}[section]
\newtheorem{lem}[thm]{Lemma}

\newtheorem{prop}[thm]{Proposition}
\theoremstyle{definition}

\theoremstyle{remark}

\newtheorem{rems}[thm]{\textbf{Remarks}}

 \makeatletter
    
    \@addtoreset{equation}{section}
  \makeatother

\makeatletter
      \def\@makefnmark{%
         \leavevmode
            \raise.9ex\hbox{\check@mathfonts
                \fontsize\sf@size\z@\normalfont%
                            \@thefnmark}%
       }
      \makeatother

\newcommand{\D}{\textrm{div}}

\newcommand{\dd}{\textrm{d}}

\begin{document}

\title[]{Liouville theorems for the Stokes equations with applications to large time estimates}
\author[]{K. Abe}
\date{}
\address[K. ABE]{Department of Mathematics, Graduate School of Science, Osaka City University, 3-3-138 Sugimoto, Sumiyoshi-ku Osaka, 558-8585, Japan}
\email{kabe@sci.osaka-cu.ac.jp}

\date{\today}

\maketitle

\begin{abstract}
We study Liouville theorems for the non-stationary Stokes equations in exterior domains in $\mathbb{R}^{n}$ under decay conditions for spatial variables. As applications, we prove that the Stokes semigroup is a bounded analytic semigroup on $L^{\infty}_{\sigma}$ of angle $\pi/2$ for $n\geq 3$. We also prove large time estimates for $n=2$ with zero net force. 
\end{abstract}

%contents
%\tableofcontents

%セクション1
\section{Introduction}

\vspace{10pt}
We consider the Stokes equations:

\begin{equation*}
\begin{aligned}
\partial_t v-\Delta{v}+\nabla{q}= 0,\quad \D\ v&=0  \qquad \textrm{in}\ \Omega\times (0,\infty),  \\
v &=0\qquad \textrm{on}\ \partial\Omega\times (0,\infty), \\
v&=v_0\hspace{18pt} \textrm{on}\ \Omega\times\{t=0\},
\end{aligned}
\tag{1.1}
\end{equation*}\\
for exterior domains $\Omega\subset \mathbb{R}^{n}$, $n\geq 2$. It is known that the Stokes semigroup $S(t):v_0\longmapsto v(\cdot,t)$ is an analytic semigroup on $L^{p}_{\sigma}$ for $p\in (1,\infty)$, of angle $\pi/2$ \cite{Sl76}, \cite{G81}, i.e., $S(t)$ is a holomorphic function in the half plane $\{ \textrm{Re}\ t>0\}$ on $L^{p}_{\sigma}$. Here, $L^{p}_{\sigma}$ denotes the $L^{p}$-closure of $C_{c,\sigma}^{\infty}$, the space of smooth solenoidal vector fields with compact support in $\Omega$. We say that an analytic semigroup is a bounded analytic semigroup of angle $\pi/2$ if the semigroup is bounded in the sector $\Sigma_{\theta}=\{t\in \mathbb{C}\backslash \{0\}\ |\ |\arg{t}|<\theta \}$ for each $\theta\in [0,\pi/2)$. The boundedness in the sector implies the bounds on the positive real line

\begin{align*}
||S(t)||\leq C,\quad ||AS(t)||\leq \frac{C'}{t}\quad t>0,  \tag{1.2}
\end{align*}\\
where $||\cdot ||$ denotes an operator norm and $A$ is the generator. The estimates (1.2) are important to study large time behavior of solutions to (1.1). In terms of the resolvent, the boundedness of $S(t)$ of angle $\pi/2$ is equivalent to the estimate 

\begin{align*}
||(\lambda-A)^{-1}||\leq \frac{C}{|\lambda|}\quad \lambda\in \Sigma_{\theta+\pi/2}.   \tag{1.3}
\end{align*}\\
When $\Omega$ is a half space, $S(t)$ is a bounded analytic semigroup on $L^{p}_{\sigma}$ of angle $\pi/2$ \cite{McCracken}, \cite{Ukai}, \cite{BM88}. The problem becomes more difficult when $\Omega$ is an exterior domain. For $n\geq 3$, the boundedness of $S(t)$ on $L^{p}_{\sigma}$ is proved in \cite{BS} based on the resolvent estimate

\begin{align*}
|\lambda|||v||_{L^{p}}+|\lambda|^{1/2} ||\nabla v||_{L^{p}}+||\nabla^{2} v||_{L^{p}}\leq C|| f||_{L^{p}}\quad 1<p<\frac{n}{2}   \tag{1.4}
\end{align*}\\
for $v=(\lambda-A)^{-1}f$ and $\lambda\in \Sigma_{\theta+\pi/2}\cup\{0\}$. The estimate (1.4) implies (1.3) for $p\in (1,n/2)$ and the case $p\in [n/2, \infty)$ follows from a duality. Due to the restriction on $p$, the two-dimensional case is more involved. Indeed, the estimate (1.4) is optimal in the sense that

\begin{align*}
||\nabla^{2} v||_{L^{p}}\leq C ||A v||_{L^{p}}\quad v\in D(A)  
\end{align*}\\
is not valid for any $p\in [n/2,\infty)$ \cite{BM92}. Here, $D(A)$ denotes the domain of $A$ on $L^{p}_{\sigma}$. For $n=2$, the boundedness of the Stokes semigroup on $L^{p}_{\sigma}$ is proved in \cite{BV} based on layer potentials for the Stokes resolvent (see also \cite{Varnhorn}).

We study the case $p=\infty$. When $\Omega$ is a half space, $S(t)$ is a bounded analytic semigroup on $L^{\infty}_{\sigma}$ of angle $\pi/2$ \cite{DHP}, \cite{Sl03}. Here, $L^{\infty}_{\sigma}$ denotes the space of all divergence-free vector fields on $L^{\infty}$, whose normal trace is vanishing on $\partial\Omega$. For bounded domains \cite{AG1} and exterior domains \cite{AG2}, analyticity of the semigroup on $L^{\infty}_{\sigma}$ follows from the a priori estimate

\begin{align*}
||v||_{L^{\infty}}+t^{1/2}||\nabla v||_{L^{\infty}}+t||\nabla^{2}v||_{L^{\infty}}+t||\partial_t v||_{L^{\infty}}+t||\nabla q||_{L^{\infty}}\leq C||v_0||_{L^{\infty}}  \tag{1.5}
\end{align*}\\
for $v=S(t)v_0$ and $t\leq T$. The estimate (1.5) implies (1.2) for $t\leq T$ and that $S(t)$ is analytic on $L^{\infty}_{\sigma}$. Moreover, the angle of analyticity is $\pi/2$ by the resolvent estimate on $L^{\infty}_{\sigma}$ \cite{AGH}. When $\Omega$ is bounded, the sup-norms in (1.5) exponentially decay as $t\to\infty$ and $S(t)$ is a bounded analytic semigroup on $L^{\infty}_{\sigma}$ of angle $\pi/2$. For exterior domains, it is non-trivial whether the Stokes semigroup is a bounded analytic semigroup on $L^{\infty}_{\sigma}$.

For the Laplace operator or general elliptic operators, it is known that corresponding semigroups are analytic on $L^{\infty}$ of angle $\pi/2$ \cite{Masuda72b}, \cite{Ste74}, \cite{Lunardi}. Moreover, if the operators are uniformly elliptic, by Gaussian upper bounds for complex time heat kernels, the semigroups are bounded analytic on $L^{\infty}$ of angle $\pi/2$; see \cite[Chapter 3]{Davies}. In particular, the heat semigroup with the Dirichlet boundary condition in an exterior domain for $n\geq 2$ is a bounded analytic semigroup on $L^{\infty}$ of angle $\pi/2$. For the Stokes equations, the Gaussian upper bound may not hold. See \cite{DHP}, \cite{Sl03} for a half space.

Large time $L^{\infty}$-estimates of the Stokes semigroup have been studied for $n\geq 3$. Maremonti \cite{Mar14} proved the estimate 

\begin{align*}
||S(t)v_0||_{L^{\infty}}\leq C||v_0||_{L^{\infty}}\quad t>0,  \tag{1.6}
\end{align*} \\
for exterior domains and $n\geq 3$ based on the short time estimate in  \cite{AG1}. Subsequently, Maremonti-Hieber \cite{MH16} proved the estimate $t||AS(t)v_0||_{L^{\infty}}\leq C||v_0||_{L^{\infty}}$ for $t>0$ and the results are extended in \cite{BH} for complex time $t\in \Sigma_{\theta}$ and $\theta\in [0,\pi/2)$ based on the approach in \cite{Mar14}. Of these papers, the case $n=2$ is excluded. We are able to observe the difference between $n\geq 3$ and $n=2$ from the representation formula of the Stokes flow due to Mizumachi \cite{Mizumachi84}; see below (1.9). In this paper, we study large time behavior of Stokes flows for $n\geq 2$ by a different approach.

Our approach is by a Liouville theorem. A Liouville theorem is a fundamental property to study regularity problems. It rules out non-trivial solutions defined in $\Omega\times (-\infty,0]$, called ancient solutions. See \cite{KNSS}, \cite{SS} for Liouville theorems of the Navier-Stokes equations and \cite{JSS} for the Stokes equations. Liouville theorems are also important to study large time behavior of solutions. In this paper, we prove non-existence of ancient solutions of (1.1) in exterior domains under spatial decay conditions. We then apply our Liouville theorems and prove the large time estimate (1.6) for complex time $t\in \Sigma_{\theta}$ and $\theta\in [0,\pi/2)$. 

We say that $v\in L^{1}_{\textrm{loc}}(\overline{\Omega}\times (-\infty,0])$ is an ancient solution to the Stokes equations (1.1) if $\D\ v=0$ in $\Omega\times (-\infty,0)$ in distributional sense and 

\begin{align*}
\int_{-\infty}^{0}\int_{\Omega}v\cdot (\partial_t\varphi+\Delta\varphi)\dd x\dd t=0,   \tag{1.7}
\end{align*}\\
for all $\varphi\in C^{2,1}_{c}(\overline{\Omega}\times (-\infty,0])$ satisfying $\D\ \varphi=0$ in $\Omega\times (-\infty,0)$ and $\varphi=0$ on $\partial\Omega\times (-\infty,0)\cup\Omega\times \{t=0\}$. Our first result is: 

\vspace{10pt}

%thm1.1
\begin{thm}[Liouville theorem]
Let $\Omega$ be an exterior domain with $C^{3}$-boundary in $\mathbb{R}^{n}$, $n\geq 2$. Let $v$ be an ancient solution to the Stokes equations (1.1). Assume that 

\begin{align*}
v\in L^{\infty}(-\infty,0; L^{p})\quad \textrm{for}\ p\in (1,\infty).  \tag{1.8}
\end{align*}\\
Then, $v\equiv 0$.
\end{thm}

\vspace{10pt}

If one removes the spatial decay condition (1.8), the assertion of Theorem 1.1 becomes false for $n\geq 3$ due to existence of stationary solutions which are asymptotically constant as $|x|\to\infty$ \cite{BM92} (cf. \cite{JSS} for bounded ancient solutions.) For $n=2$, it is known that bounded stationary solutions do not exist \cite{ChangFinn}. We show that ancient solutions in $L^{\infty}(-\infty,0; L^{p})$ are extendable to bounded entire functions by using analyticity of the Stokes semigroup. The Liouville theorem then follows from the information of the kernel of the Stokes operator on $L^{p}_{\sigma}$.

Theorem 1.1 is useful to study the large time estimate (1.6) for $t>0$. We invoke the representation formula of the Stokes flow 

\begin{align*}
v(x,t)=\int_{\Omega}\Gamma (x-y,t)v_0(y)\dd y+\int_{0}^{t}\int_{\partial\Omega}V(x-y,t-s)T(y,s)n(y)\dd H(y)\dd s.  \tag{1.9}
\end{align*}\\
Here, $T=\nabla v+\nabla^{T} v-qI$ is the stress tensor and $V=(V_{ij}(x,t))_{i,j}$ is the Oseen tensor

\begin{align*}
V_{ij}(x,t)=\Gamma(x,t)\delta_{ij}-\partial_i\partial_j\int_{\mathbb{R}^{n}}E(x-y)\Gamma(y,t)\dd y,   \tag{1.10}
\end{align*}\\
defined by the heat kernel $\Gamma$ and the fundamental solutions of the Laplace equation $E$. The formula (1.9) is obtained by regarding $v=S(t)v_0$ as the Stokes flow in $\mathbb{R}^{n}$ with a measure as the external force; see Remarks 3.5 (i). It describes the asymptotic behavior of bounded Stokes flow as $|x|\to\infty$. We show that if the Stokes flow is bounded for all $t>0$, the stress tensor is also bounded on $\partial\Omega$. Observe that by the pointwise estimate of the Oseen tensor

\begin{align*}
|V(x,t)|\leq \frac{C}{(|x|+t^{1/2})^{n}},\quad x\in \mathbb{R}^{n},\ t>0,   \tag{1.11}
\end{align*}\\
the remainder term is estimated by 

\begin{align*}
\left|v(x,t)-\int_{\Omega}\Gamma(x-y,t)v_0(y)\dd y\right|
\leq \frac{C}{|x|^{n-2}}\sup_{0<s\leq t}||T||_{L^{\infty}(\partial\Omega)}(s),  \tag{1.12}
\end{align*}\\
for $|x|\geq 2R_0$ and $t>0$ with $R_0\geq \textrm{diam}\ \Omega^{c}$. The right-hand side is decaying as $|x|\to\infty$ uniformly for $t>0$ if $n\geq 3$. We show that the large time estimate (1.6) is reduced to showing non-existence of ancient solutions by a contradiction argument. Since the remainder term estimate (1.12) yields a decay condition for ancient solutions as $|x|\to\infty$, we are able to obtain a contradiction by applying the Liouville theorem (Theorem 1.1).

We further extend our approach to obtain (1.6) for complex time $t\in \Sigma_{\theta}$ and $\theta\in [0,\pi/2)$. To this end, we consider ancient solutions in the sector $\Lambda=\{t\in \mathbb{C}\backslash\{0\}\ |\ \pi-\theta\leq \arg{t}\leq -\pi/2 \}$. We set the segment $I_{T}=\{t\in \Lambda\ |\ -T\leq \textrm{Re}t\leq 0,\ \textrm{Im}\ t=-T\tan{\theta}  \}$ for $T>0$. We say that $v$ is an ancient solution in $\Omega\times \Lambda$ if $v$ satisfies the Stokes equations (1.1) on each segment $I_T$ for $T>0$ in a weak sense, i.e., $v\in L^{1}_{\textrm{loc}}(\overline{\Omega}\times I_{T})$ satisfies $\D\ v=0$ in $\Omega\times I_{T}$ in the distributional sense, 

\begin{align*}
\int_{-T}^{0}\hspace{-3pt}\int_{\Omega}v(x,\alpha+i\beta)(\partial_\alpha\varphi+\Delta \varphi)\dd x\dd \alpha=-\int_{\Omega}v(x,-T+i\beta)\varphi(x,-T)\dd x,\hspace{3pt}  \beta=-T\tan{\theta}   \tag{1.13}
\end{align*}\\
for all $\varphi\in C^{2,1}_{c}(\overline{\Omega}\times [-T,0])$ satisfying $\D\ \varphi=0$ in $\Omega\times (-T,0)$ and $\varphi=0$ on $\partial\Omega\times (-T,0)\cup \Omega\times \{t=0\}$. For $\theta=0$, an  ancient solution in the sector is an ancient solution on the negative line. We prove non-existence of ancient solutions in the sector under the condition $v\in L^{\infty}(\Lambda; L^{p})$ for $p\in (1,\infty)$ (Theorem 2.7). Since the formula (1.9) is extendable for complex time, we apply the Liouville theorem in the sector and obtian (1.6) for $t\in \Sigma_{\theta}$ and $\theta\in [0,\pi/2)$. We now state our main results.

\vspace{10pt}

%thm1.2
\begin{thm}
When $n\geq 3$, the Stokes semigroup is a bounded analytic semigroup on $L^{\infty}_{\sigma}$ of angle $\pi/2$. 
\end{thm}

\vspace{10pt}
For $n=2$, the remainder term estimate (1.12) is different. By a simple calculation from the formula (1.9), we see an asymptotic profile of the two-dimensional Stokes flow:

\begin{align*}
\left|v(x,t)-\int_{\Omega}\Gamma(x-y,t)v_0(y)\dd y-\int_{0}^{t}V(x,t-s)N(s)\dd s\right|\leq\frac{C}{|x|} \sup_{0<s\leq t}||T||_{L^{\infty}(\partial\Omega)}(s),  \tag{1.14}
\end{align*}\\
for $|x|\geq 2R_0$ and $t>0$, with the net force

\begin{align*}
N(s)=\int_{\partial\Omega}T(y,s)n(y)\dd H(y).
\end{align*}\\
Since $|\int_{0}^{t}V(x,s)\dd s|\lesssim \log{(1+t/|x|^{2})}$, the decay as $|x| \to\infty$ of the third term in the left-hand side is not uniform for $t>0$ in contrast to (1.12) for $n\geq 3$. If the net force vanishes, the situation is the same as $n=3$ and we are able to prove (1.6). For example, if initial data is rotationally symmetric, the net force vanishes. Following \cite{Bran}, we consider initial data invariant under a cyclic group or a dihedral group. For integers  $m\geq 2$, we set the matrices

\begin{align*}
R_m=
\begin{pmatrix}
\cos(2\pi/m) &  -\sin{(2\pi/m)} \\
\sin{(2\pi/m)} & \cos{(2\pi/m)}	 
\end{pmatrix},
\quad
J=
\begin{pmatrix}
1 &  0 \\
0 & -1	 
\end{pmatrix}.
\end{align*}\\
Let $C_m$ denote the cyclic group of order $m$ generated by the rotation $R_m$. Let $D_m$ denote the dihedral group of order $2m$ generated by $R_m$ and the reflection $J$. Any finite subgroup of the orthogonal group $O(2)$ is either a cyclic group or a dihedral group. See \cite[Chapter 2]{GroveBenson}. Let $G$ be a subgroup of $O(2)$ and $\Omega^{c}$ be a disk centered at the origin. We say that a vector field $v$ is $G$-covariant if $v(x)={}^{t}Av(Ax)$ for all $A\in G$ and $x\in \Omega^{c}$. It is known that if $v_0$ is $C_m$-covariant, so is $v=S(t)v_0$ and the net force vanishes \cite{HeMiyakawa06}. Thus for $C_m$-covariant vector fields $v_0 \in L^{\infty}_{\sigma}$, the remainder term estimate is the same as $n=3$.

\vspace{10pt}

\begin{thm}
For $n=2$, the estimate (1.6) holds for $t\in \Sigma_{\theta}$ and $v_0\in L^{\infty}_{\sigma}$, for which the net force vanishes (e.g., $C_m$-covariant vector fields when $\Omega^{c}$ is a disk.)
\end{thm}

\vspace{10pt}
Theorem 1.3 improves the pointwise estimates of the two-dimensional Navier-Stokes flows for rotationally symmetric initial data \cite{HeMiyakawa06}, in which the estimate (1.6) is noted as an open question together with the applications to the nonlinear problem. We are able to apply (1.6) to improve the results although initial data is restricted to rotationally symmetric; see Remarks 5.8 (iii). 

We hope it is possible to extend our approach to study the case with net force, for which (1.6) is unknown even if initial data is with finite Dirichlet integral. The estimate (1.6) with net force is important to study large time behavior of asymptotically constant solutions as $|x|\to\infty$. We refer to \cite{A7} for asymptotically constant solutions of the two-dimensional Navier-Stokes equations. See also \cite{MaremontiShimizu}.\\

This paper is organized as follows. In Section 2, we prove Theorem 1.1. In Section 3, we prove the remainder term estimate (1.12). In Section 4, we prove (1.6) for positive time. In Section 5, we prove (1.6) for complex time and complete the proof of Theorems 1.2. After the proof of Theorem 1.2, we prove Theorem 1.3.

\vspace{15pt}

\section{Liouville theorems on $L^{p}$}

\vspace{10pt}

We prove Theorem 1.1. We show that ancient solutions are bounded entire functions on $L^{p}$ under the condition (1.8). To this end, we prove a uniqueness theorem and extend ancient solutions $v(t)$ for $\textrm{Re}\ t>-T$ for each $T>0$ by the analytic continuation $v(t)=S(t+T)v(-T)$. After the proof of Theorem 1.1, we prove a Liouville theorem in a sector (Theorem 2.7).

\vspace{10pt}

%lem2.1
\begin{lem}[Uniqueness]
Let $v\in L^{1}_{\textrm{loc}}(\overline{\Omega}\times [0,T])$ satisfy $\D\ v=0$ and 

\begin{align*}
\int_{0}^{T}\int_{\Omega} v\cdot (\partial_t \varphi+\Delta \varphi)\dd x\dd t=0  \tag{2.1}
\end{align*}\\
for all $\varphi\in C^{2,1}_{c}(\overline{\Omega}\times [0,T])$ such that $\D\ \varphi=0$ in $\Omega\times (0,T)$, $\varphi=0$ on $\partial\Omega\times (0,T)$ and $\Omega\times \{t=T\}$. Assume that 

\begin{align*}
v\in L^{\infty}(0,T; L^{p})\quad \textrm{for}\ p\in (1,\infty).   \tag{2.2}
\end{align*}\\
Then, $v\equiv 0$.
\end{lem}

\vspace{10pt}
We first extend test functions of (2.1) under the boundedness (2.2). We use the Bogovski\u\i\ operator.
\vspace{10pt}

%prop2.2
\begin{prop}
(i) Let $D=\{x\in \mathbb{R}^{n}\ |\ 1<|x|<2\}$ and $L^{p}_{\textrm{av}}(D)=\{h\in L^{p}(D)\ |\ \int_{D}h\dd x=0\}$. There exists a bounded linear operator $B:L^{p}_{av}(D)\longrightarrow W^{1,p}_{0}(D)$, $p\in (1,\infty)$, such that $w=B(h)$ satisfies 

\begin{align*}
\D\ w=h\quad \textrm{in}\ D,\quad 
w =0\quad \textrm{on}\ \partial D.  \tag{2.3}
\end{align*}\\
Moreover, the operator $B$ acts as a bounded operator from $W^{k,p}_{0}(D)$ to $W^{k+1,p}_{0}(D)$ for positive integers $k$.

\noindent 
(ii) Let $D_R=\{R<|x|<2R  \}$. There exists a bounded operator $B_R: L^{p}_{\textrm{av}}(D_R)\longrightarrow W^{1,p}_{0}(D_R)$ satisfying (2.3) in $D_R$. Moreover, the estimate

\begin{align*}
||\nabla^{k+1}B_R(h)||_{L^{p}(D_R)}\leq C||\nabla^{k}h||_{L^{p}(D_R)} \tag{2.4}
\end{align*}\\
holds with some constant $C$, independent of $R>0$.
\end{prop}

\vspace{5pt}

\begin{proof}
See \cite{Bogovskii79}, \cite{BS90}, \cite[Theorem III.3.3]{Gal} for the assertion (i). The operator $B_R$ is constructed by (i) and dilation.   
\end{proof}

\vspace{10pt}

%prop2.3
\begin{prop}
Under the assumption of Lemma 2.1, the equality (2.1) is extendable for all $\varphi\in C^{2,1}(\overline{\Omega}\times [0,T])$ such that $\D\ \varphi=0$ in $\Omega\times (0,T)$, $\varphi=0$ on $\partial\Omega\times (0,T)$ and $\Omega\times \{t=T\}$, 

\begin{align*}
\partial_t^{s}\partial_x^{k}\varphi \in L^{\infty}(0,T; L^{q})   \tag{2.5}
\end{align*}\\
for $2s+|k|\leq 2$ and $1/p+1/q=1$.
\end{prop}

\vspace{10pt}

\begin{proof}
Let $B_0(R_0)$ denote an open ball centered at the origin with radius $R_0>0$. We take $R_0>0$ so that $\Omega^{c}\subset B_{0}(R_0)$. Let $\theta\in C^{\infty}_{c}[0,\infty)$ be a function such that $\theta\equiv 1$ in $[0,1]$ and $\theta\equiv 0$ in $[2,\infty)$. We set $\theta_{R}(x)=\theta(|x|/R)$ for $R\geq R_0$. Since $\D\ \varphi=0$ in $\Omega$ and $\varphi=0$ on $\partial\Omega$, the average of $h_R=\varphi\cdot \nabla \theta_R$ in $D_R$ is zero. We set $w_R=B_R(h_R)$ and consider its zero extension to $\mathbb{R}^{n}$. Since $\varphi\in C^{2,1}(\overline{\Omega}\times [0,T])$ and $B_R$ is a  linear operator, we see that $w_R\in C^{2,1}_{c}(\mathbb{R}^{n}\times [0,T])$ by (2.4) and the Sobolev embedding. We set
 
\begin{align*}
\varphi_{R}=\varphi\theta_{R}-w_R
\end{align*}\\
so that $\varphi_R\in C^{2,1}_{c}(\overline{\Omega}\times [0,T])$ satisfies $\D\ \varphi_R=0$ in $\Omega$ and $\varphi_R=0$ on $\partial\Omega$ for $t\in [0,T]$. By substituting $\varphi_R$ into (2.1), we see that 

\begin{align*}
0=\int_{0}^{T}\int_{\Omega}v\cdot (\partial_t \varphi\theta_R+\Delta\varphi\theta_R+2\nabla \varphi\cdot \nabla \theta_R+\varphi\Delta\theta_R)\dd x\dd t
-\int_{0}^{T}\int_{\Omega}v\cdot (\partial_t w_R+\Delta w_R)\dd x\dd t.
\end{align*}\\
By (2.2) and (2.5), the first term converges to the integral of $v\cdot (\partial_t\varphi+\Delta \varphi)$ in $\Omega\times (0,T)$ as $R\to\infty$. We show that the second term converges to zero. We show the convergence of the integral of $v\cdot \partial_t w_R$ since $v\cdot \Delta w_R$ is estimated by a similar way. Since $\partial_t w_R=B(\partial_t h_R)$, by the Poincar\'e inequality \cite{E} and (2.4) we estimate

\begin{align*}
||\partial_t w_R||_{L^{q}(D_R)}
\leq CR||\nabla \partial_t w_R||_{L^{q}(D_R)}
&=CR||\nabla B_R(\partial_t h_R)||_{L^{q}(D_R)} \\
&\leq C'R||\partial_t h_R||_{L^{q}(D_R)}\\
&\leq C''||\partial_t \varphi||_{L^{q}(D_R)}.
\end{align*}\\
It follows that

\begin{align*}
\left|\int_{0}^{T}\int_{\Omega}v\cdot \partial_t w_R \dd x\dd t\right|
\leq C||v||_{L^{\infty}(0,T;L^{p})} \left(\int_{0}^{T}||\partial_t\varphi||_{L^{q}(D_R)}\dd t\right)\to0\quad R\to\infty.
\end{align*}\\
The proof is complete.
\end{proof}

\vspace{10pt}
We apply a duality argument to prove the uniqueness. To this end, we show existence of solutions to the adjoint problem.
\vspace{10pt}

%prop2.4
\begin{prop}
For $f\in C^{\infty}_{c}(\Omega\times (0,T))$ satisfying $\D\ f=0$, there exists a solution $(\varphi,\nabla \pi)\in C^{2,1}(\overline{\Omega}\times [0,T])\times C(\overline{\Omega}\times [0,T])$ of 

\begin{equation*}
\begin{aligned}
\partial_t \varphi+\Delta\varphi-\nabla \pi=f,\quad \D\ \varphi=0\quad \textrm{in}\ \Omega\times (0,T),\\
\varphi=0\quad \textrm{on}\ \partial\Omega\times(0,T)\cup \Omega\times \{t=T\},
\end{aligned}
\tag{2.6}
\end{equation*}\\
satisfying (2.5) and $\nabla \pi\in L^{\infty}(0,T; L^{q})$ for all $q\in (1,\infty)$.
\end{prop}

\vspace{10pt}

%prop2.5
\begin{prop}
For $g\in C^{\infty}_{c}(\Omega\times (0,T))$ satisfying $\D\ g=0$, there exists a solution $(\psi,\nabla s)\in C^{2,1}(\overline{\Omega}\times [0,T])\times C(\overline{\Omega}\times [0,T])$ of 

\begin{equation*}
\begin{aligned}
\partial_t \psi-\Delta\psi+\nabla s=g,\quad \D\ \varphi=0\quad \textrm{in}\ \Omega\times (0,T),\\
\psi=0\quad \textrm{on}\ \partial\Omega\times(0,T)\cup \Omega\times \{t=0\},
\end{aligned}
\tag{2.7}
\end{equation*}\\
such that $\partial_t^{s}\partial_x^{k}\psi, \nabla s\in L^{\infty}(0,T; L^{q})$ for $2s+|k|\leq 2$ and all $q\in (1,\infty)$.
\end{prop}

\vspace{5pt}

\begin{proof}
Let $\mathbb{P}$ denote the Helmholtz projection operator on $L^{p}$ \cite{SiSo}. By the Stokes operator $A=\mathbb{P}\Delta$ with the domain $D(A)=W^{2,p}\cap W^{1,p}_{0}\cap L^{p}_{\sigma}$ and the Stokes semigroup $S(t)=e^{tA}$, we set 
\begin{align*}
\psi(x,t)=\int_{0}^{t}S(t-s)g(s)\dd s.   
\end{align*}\\
Since $g$ is smooth and $D(A)\subset W^{2,p}$, we see that $\partial_t^{s}\partial_x^{k}\psi\in L^{\infty}(0,T; L^{p})$ for $2s+|k|\leq 2$. We set $\nabla s=(1-\mathbb{P})\Delta \psi$. Since $\partial_{t}\psi-A\psi=g$ on $L^{p}$, $(\psi,\pi)$ satisfies (2.7).

It remains to show that $(\psi,s)$ is continuous up to second orders in $\overline{\Omega}\times [0,T]$. Since $g$ is smooth, in particular $\psi\in C^{1}([0,T]; D(A^{2}))$. We take bounded domains $\Omega''$ and $\Omega'$ such that $\Omega''\subset \Omega'\subset \Omega$. Since the boundary is $C^{3}$, we apply the higher regularity estimate for the Stokes operator \cite{Gal} to estimate

\begin{align*}
||\psi||_{W^{3,p}(\Omega'')}+||s||_{W^{2,p}(\Omega'')}
\leq C(||A\psi||_{W^{1,p}(\Omega')}+||\psi||_{W^{1,p}(\Omega')}+||s||_{L^{p}(\Omega')} ).
\end{align*}\\
Thus, $\psi\in C([0,T]; W^{3,p}_{\textrm{loc}}(\overline{\Omega}))$. By the Sobolev embedding for $p>n$, we see that $\nabla^{2}\psi$ is continuous in $\overline{\Omega}\times [0,T]$. Since $\partial_t \psi\in C([0,T]; D(A^{2}))$ and by (2.7), $\partial_t\psi$ and $\nabla s$ are continuous in $\overline{\Omega}\times [0,T]$. The proof is complete.
\end{proof}

\vspace{10pt}

\begin{proof}[Proof of Proposition 2.4]
For $g(x,t)=-f(x,T-t)$, we take $(\psi,s)$ by Proposition 2.5 and set $\varphi(x,t)=\psi(x,T-t)$, $\pi(x,t)=s(x,T-t)$. Then, $(\varphi,\pi)$ satisfies the desired properties.
\end{proof}

\vspace{10pt}

\begin{proof}[Proof of Lemma 2.1]
For an arbitrary $f\in C^{\infty}_{c}(\Omega\times (0,T))$ satisfying $\D\ f=0$, we take a solution $(\varphi,\pi)$ of the adjoint problem (2.6) by Proposition 2.4. Since we extended test functions in Proposition 2.3, it follows that

\begin{align*}
\int_{0}^{T}\int_{\Omega}v\cdot f \dd x\dd t=\int_{0}^{T}\int_{\Omega}v\cdot (\partial_t\varphi+\Delta \varphi-\nabla \pi)\dd x\dd t
=0.
\end{align*}\\
Since $f$ is solenoidal, by deRham's theorem \cite{deRham}, \cite[Theorem 1.1]{SiSo}, there exists a function $\Phi$ such that $v=\nabla \Phi$. Since $\nabla \Phi$ is harmonic and vanishing on $\partial\Omega$ and $|x|\to\infty$, applying the maximum principle implies $\nabla \Phi\equiv 0$. The proof is complete. 
\end{proof}

\vspace{10pt}

%prop2.6
\begin{prop}
The kernel of the Stokes operator on $L^{p}_{\sigma}$ for $p\in (1,\infty)$ is zero, i.e., $N(A)=\{v\in D(A)\ |\ Av=0\}=\{0\}$. 
\end{prop}

\vspace{10pt}
\begin{proof}
See \cite[Corollary 3.6]{GS} for $n\geq 3$ and \cite{ChangFinn}, \cite[p.297]{BV} for $n=2$.
\end{proof}

\vspace{10pt}

\begin{proof}[Proof of Theorem 1.1]
We take an arbitrary $T>0$ and set $u(t)=v(t-T)$ for $t\in [0,T]$ so that $u$ satisfies 

\begin{align*}
\int_{0}^{T}\int_{\Omega}u\cdot (\partial_t \varphi+\Delta \varphi)\dd x\dd t=-\int_{\Omega}u(x,0)\cdot \varphi(x,0)\dd x  \tag{2.8}
\end{align*}\\
for all $\varphi\in C^{2,1}_{c}(\overline{\Omega}\times [0,T])$, $\D\ \varphi=0$ in $\Omega\times (0,T)$, $\varphi=0$ on $\partial\Omega\times (0,T)\cup \Omega\times \{t=T\}$. We set $\tilde{u}(t)=S(t)u(0)$. Since $S(t)$ is a bounded analytic semigroup of angle $\pi/2$ on $L^{p}_{\sigma}$, $\tilde{u}(t)$ is defined for $\textrm{Re}\ t>0$ and bounded in $\Sigma_{\theta}$ for $\theta\in [0,\pi/2)$. Since $\tilde{u}(t)=S(t)u(0)$ also satisfies (2.8), applying Lemma 2.1 implies $u(t)=S(t)u(0)$. Thus $u(t)$ is continued for $\textrm{Re}\ t>0$ and satisfies 

\begin{align*}
\sup\left\{||u||_{L^{p}}(t)\ \middle|\  t\neq 0,\ |\arg{t}|\leq \theta \right\}\leq C||u||_{L^{p}}(0).
\end{align*}\\
This means that 

\begin{align*}
\sup\left\{||v||_{L^{p}}(t)\  \middle|\ t\neq -T,\ \arg{(t+T)}|\leq \theta \right\}\leq C||v||_{L^{p}}(-T).
\end{align*}\\
Since the right-hand side is uniformly bounded for $T>0$ by (1.8), the ancient solution is a bounded entire function. Thus, $\partial_t v\equiv 0$ by the Liouville theorem and $v\equiv 0$ follows from Proposition 2.6.
\end{proof}

\vspace{10pt}
We extend Theorem 1.1 to a sector.
\vspace{10pt}

%thm 2.7
\begin{thm}
Let $v$ be an ancient solution to (1.1) in $\Omega\times \Lambda$ for $\theta\in [0,\pi/2)$. Assume that 

\begin{align*}
v\in L^{\infty}(\Lambda; L^{p})\quad  \textrm{for}\ p\in (1,\infty). \tag{2.9}
\end{align*}\\
Then, $v\equiv 0$.
\end{thm}

\vspace{10pt}

\begin{proof}
It suffices to show the case $\theta\neq 0$. We take an arbitrary $T>0$ and set $\tilde{T}=T+iT\tan{\theta}$. By translation, we set $u(t)=v(t-\tilde{T})$. Then by (1.13), $u$ satisfies 

\begin{align*}
\int_{0}^{T}\int_{\Omega}u(x,\alpha)(\partial_\alpha\varphi+\Delta \varphi)\dd x\dd \alpha=-\int_{\Omega}u(x,0)\varphi(x,0)\dd x, 
\end{align*}\\
for $\varphi\in C^{2,1}_{c}(\overline{\Omega}\times [0,T))$ satisfying $\D\ \varphi=0$ in $\Omega\times (0,T)$ and $\varphi=0$ on $\partial\Omega\times (0,T)\cup \Omega\times \{t=T\}$. The by the uniqueness theorem as in the proof of Theorem 1.1, $u(t)=S(t)u(0)$ is continued to $\textrm{Re}\ t>0$ and bounded in $\Sigma_{\theta'}$ for $\theta'\in [0,\pi/2)$. We take $\theta'\in (\theta,\pi/2)$ and estimate 

\begin{align*}
\sup\left\{||u||_{L^{p}}(t) \middle|\ t\neq 0,\ |\arg{t}|\leq \theta' \right\}\leq C||u||_{L^{p}}(0). 
\end{align*}\\
Since $v(t)=S(t+\tilde{T})v(-\tilde{T})$, $v$ is continued to $\{\textrm{Re}\ t>-T\}$ and 

\begin{align*}
\sup\left\{||v||_{L^{p}}(t) \middle|\ t\neq -\tilde{T},\ |\arg({t+\tilde{T}})|\leq \theta' \right\}\leq C||v||_{L^{p}}(-\tilde{T}).  
\end{align*}\\
Since the right-hand side is uniformly bounded for $T>0$ by (2.9), $v$ is a bounded entire function. Thus $v\equiv 0$ follows.
\end{proof}

\vspace{15pt}

\section{Representation formula}

\vspace{10pt}

We prove the remainder term estimate (1.12) from the representation formula (1.9). We also prove an $L^{\infty}$-estimate of the stress tensor (3.6). The estimates (1.12) and (3.6) are used in the next section in order to prove (1.6) for $t>0$. We begin with the Stokes equations for bounded data.

\vspace{10pt}

\begin{prop}
(i) For $v_0\in L^{\infty}_{\sigma}$, there exists a unique solution $(v,\nabla q)\in C^{2,1}(\overline{\Omega}\times (0,T])\times C(\overline{\Omega}\times (0,T])$ of the Stokes equations (1.1) such that $v(\cdot,t)\to v_0$ a.e. in $\Omega$ as $t\to0$ and 

\begin{align*}
\sup_{0<t\leq T}\left(t^{\frac{|k|}{2}+s} ||\partial_t^{s}\partial_x^{k}v||_{L^{\infty}}(t)+t||\nabla q||_{\infty}(t)\right)\leq C||v_0||_{L^{\infty}}  \tag{3.1}
\end{align*}\\
for $T>0$ and $2s+|k|\leq 2$, with some constant $C=C(T)$, independent of $v_0$.

\noindent 
(ii) The Stokes semigroup $S(t): v_0\longmapsto v(\cdot,t)$ is an analytic semigroup on $L^{\infty}_{\sigma}$ of angle $\pi/2$ and satisfies 

\begin{align*}
\sup_{0<t\leq T}t^{\frac{|k|}{2}+s} ||\partial_t^{s}\partial_x^{k}S(t)v_0||_{L^{\infty}}\leq C||v_0||_{L^{\infty}}.   \tag{3.2}
\end{align*}\\
The semigroup $S(t)$ is weakly-star continuous on $L^{\infty}$ at time zero. 
\end{prop}

\vspace{5pt}
\begin{proof}
The assertion is proved in \cite{AG2} except the angle of analyticity \cite{AGH}.
\end{proof}

\vspace{10pt}

%lem3.2
\begin{lem}[Representation formula]
\noindent 
(i) For $v=S(t)v_0$ and $v_0\in L^{\infty}_{\sigma}$, 

\begin{equation*}
\begin{aligned}
v(x,t)=\int_{\Omega}\Gamma (x-y,t)v_0(y)\dd y+\int_{0}^{t}\int_{\partial\Omega}V(x-y,t-s)T(y,s)n(y)\dd H(y)\dd s,
\end{aligned}
\tag{3.3}
\end{equation*}\\
holds for $x\in \Omega$ and $t>0$.

\noindent
(ii) For $n\geq 3$, there exists a constant $C$ such that 

\begin{align*}
|v(x,t)|\leq ||v_0||_{L^{\infty}}+\frac{C}{d(x)^{n-2}}\sup_{0<s\leq t}||T||_{L^{\infty}(\partial\Omega)}(s), \quad x\in \Omega,\ t>0.  \tag{3.4}
\end{align*}\\
holds for $x\in \Omega$ and $t>0$, where $d(x)=\textrm{dist}(x,\partial\Omega)$ denotes the distance from the boundary.
\end{lem}

\vspace{5pt}

\begin{proof}
The formula (3.3) holds for $v_0\in C^{\infty}_{c,\sigma}$ \cite{Mizumachi84} (see also \cite[Appendix B]{HeMiyakawa06}). Since the solution $S(t)v_0$ for $v_0\in L^{\infty}_{\sigma}$ is constructed as a limit of a sequence $v_m=S(t)v_{0,m}$ for $\{v_{0,m}\}\subset  C^{\infty}_{c,\sigma}$ such that $v_{0,m}\to v_0$ a.e. in $\Omega$ \cite{AG2}, (3.3) is extendable for $v_0\in L^{\infty}_{\sigma}$. 

By the pointwise estimate of the Oseen tensor (1.11), we estimate

\begin{align*}
\left|\int_{0}^{t}\int_{\partial\Omega}V(x-y,t-s)T(y,s)n(y)\dd H(y)\dd s\right|
&\leq C\sup_{0<s\leq t}||T||_{L^{\infty}(\partial\Omega)}(s)
\int_{0}^{t}\int_{\partial\Omega}\frac{\dd H(y)\dd s}{(|x-y|+s^{1/2})^{n}}\\
&\leq \frac{C'}{d(x)^{n-2}}\sup_{0<s \leq t}||T||_{L^{\infty}(\partial\Omega)}(s)\int_{0}^{\infty}\frac{\dd s}{(1+s^{1/2})^{n}}.
\end{align*}\\
Since the right-hand side is finite for $n\geq 3$, (3.4) holds.
\end{proof}

\vspace{10pt}
We estimate the stress tensor $T$ on $\partial\Omega$ by using a bound for $v=S(t)v_0$.

\vspace{10pt}

%Prop3.3
\begin{lem}
Let $R_0>0$ be a constant such that $\Omega^{c}\subset B_0(R_0)$ and $\Omega_0=\Omega\cap B_0(R_0)$. Let $q$ be an associated pressure for $v=S(s)v_0$ and $v_0\in D(A)$ satisfying 

\begin{align*}
\int_{\Omega_0}q(x,t)\dd x=0.  \tag{3.5}
\end{align*}\\
Then, the estimate

\begin{align*}
\sup_{0<s\leq t}||T||_{L^{\infty}(\partial\Omega)}(s)\leq C\left(||v_0||_{D(A)}+\sup_{0<s\leq t}||v||_{L^{\infty}(\Omega)}(s) \right)  \tag{3.6}
\end{align*}\\
holds for all $t>0$ with some constant $C$, independent of $v_0$, where $A$ denotes the Stokes operator on $L^{\infty}_{\sigma}$ with the domain $D(A)$, equipped with the graph-norm $||\cdot ||_{D(A)}$.
\end{lem}

\vspace{10pt}
In order to prove Lemma 3.3, we use the Stokes resolvent estimate on $L^{\infty}_{\sigma}$.
\vspace{10pt}

%prop3.4
\begin{prop}
(i) For $f\in L^{\infty}_{\sigma}$, there exists a unique solution $(v,\nabla q)\in W^{2,p}_{\textrm{ul}}(\overline{\Omega})\times L^{p}_{\textrm{ul}}(\overline{\Omega})\cap L^{\infty}_{d}(\Omega)$, $p>n$, of the Stokes equations

\begin{equation*}
\begin{aligned}
v-\Delta v+\nabla q=f,\quad \D\ v&=0\quad \textrm{in}\ \Omega, \\
v&=0\quad \textrm{on}\ \partial\Omega,
\end{aligned}
\tag{3.7}
\end{equation*}\\
satisfying the estimate

\begin{align*}
||v||_{W^{2,p}_{\textrm{ul}}(\overline{\Omega}) }+||\nabla q||_{L^{p}_{\textrm{ul}}(\overline{\Omega}) }
\leq C||f||_{L^{\infty}(\Omega)}.\tag{3.8}
\end{align*}\\
Here, $L^{p}_{\textrm{ul}}(\overline{\Omega})$ denotes the uniformly local $L^{p}$ space equipped with the norm 

\begin{align*}
||g||_{L^{p}_{\textrm{ul}}(\overline{\Omega}) }
=\sup\left\{||g||_{L^{p}(B_{x_0}(1)\cap \Omega)}\ \middle|\ x_0\in \Omega\     \right\}.
\end{align*}\\
The space $W^{2,p}_{\textrm{ul}}(\overline{\Omega})$ is equipped with the norm $||v||_{W^{2,p}_{\textrm{ul}}(\overline{\Omega}) }=\sum_{|k|\leq 2}||\partial_x^{k}v||_{L^{p}_{\textrm{ul}}(\overline{\Omega}) }$ and $L^{\infty}_{d}(\Omega)$ denotes the space of functions $g\in L^{1}_{\textrm{loc}}(\Omega)$ such that $dg\in L^{\infty}(\Omega)$ for $d(x)=\textrm{dist}(x,\partial\Omega)$.

\noindent 
(ii) There exists a constant $C$ such that 

\begin{align*}
||T||_{L^{\infty}(\partial\Omega)}\leq C||f||_{L^{\infty}(\Omega)}  \tag{3.9}
\end{align*}\\
holds for solutions of (3.7) satisfying (3.5).
\end{prop}

\vspace{10pt}
\begin{proof}
The assertion (i) is proved in \cite[Theorem 1.1]{AGH}. We prove (ii). Since the average of $q$ in $\Omega_0$ is zero, by the Poincar\'e inequality \cite{E}, we estimate 

\begin{align*}
||q||_{L^{p}(\Omega_0)}
\leq C||\nabla q||_{L^{p}(\Omega_0)}
\leq C'||\nabla q||_{L^{p}_{\textrm{ul}}(\overline{\Omega})}.
\end{align*}\\
By the Sobolev inequality for $p>n$,, it follows that

\begin{align*}
||q||_{L^{\infty}(\Omega_0)}\leq C||\nabla q||_{L^{p}_{\textrm{ul}}(\overline{\Omega})}.
\end{align*}\\
Since $W^{2,p}_{\textrm{ul}}(\overline{\Omega})\subset W^{1,\infty}(\Omega)$ with continuous injection, (3.9) follows from (3.8).
\end{proof}

\vspace{10pt}

\begin{proof}[Proof of Lemma 3.3]
Since $v=S(s)v_0$ satisfies (3.7) for $f=-Av+v$, we estimate 

\begin{align*}
||T||_{L^{\infty}(\partial\Omega)}(s)\leq C||v||_{D(A)}(s),
\end{align*}\\
by (3.9). We may assume that $t>1$. For $s\in (0,1)$, it follows from (3.2) that 

\begin{align*}
||v||_{D(A)}(s)=||S(s)v_0||_{L^{\infty}}+||AS(s)v_0||_{L^{\infty}}
\leq C(||v_0||_{L^{\infty}}+||Av_0||_{L^{\infty}} )=C||v_0||_{D(A)}.
\end{align*}\\
For $s\in [1,t]$, we estimate 

\begin{align*}
||v||_{D(A)}(s)=||S(s)v_0||_{L^{\infty}}+||AS(1)S(s-1)v_0||_{L^{\infty}}
&\leq ||S(s)v_0||_{L^{\infty}}+C||S(s-1)v_0||_{L^{\infty}}  \\
&\leq C'\sup_{0<s\leq t}||v||_{L^{\infty}}(s).
\end{align*}\\
We obtained (3.6).
\end{proof}

\vspace{10pt}

%rem 3.5
\begin{rems}
\noindent
(i) The formula (3.3) is obtained by regarding $v=S(t)v_0$ as a solution of

\begin{equation*}
\begin{aligned}
\partial_tv-\Delta v+\nabla q=f\quad \D\ v&=0\quad\textrm{in}\ \mathbb{R}^{n}\times (0,\infty),\\
v&=v_0\hspace{7pt} \textrm{on}\ \mathbb{R}^{n}\times \{t=0\},
\end{aligned}
\tag{3.10}
\end{equation*}\\
for $f=Tn\chi$ and $\chi\in C_{0}(\mathbb{R}^{n})^{*}$ such that 

\begin{align*}
<\chi,\varphi>=\int_{\partial\Omega}\varphi(y)\dd {H}(y).
\end{align*}\\
To see this, we use a cut-off function $\eta_{\delta}(x)=\eta(d(x)/\delta)$ for $\delta>0$ by a smooth function $\eta$ such that $\eta\equiv 0$ in $[0,1]$ and $\eta\equiv 1$ in $[2,\infty)$. We set $v_{\delta}=v\eta_{\delta}+\nabla \Phi_{\delta}$ and $q_{\delta}=q\eta_{\delta}-(\partial_t-\Delta) \Phi_{\delta}$ by $\Phi_{\delta}=E*(v\cdot \nabla\eta_{\delta})$ so that $(v_\delta,q_{\delta})$ satisfies (3.10) for $f_{\delta}=-\nabla v\nabla \eta_{\delta}+q\nabla \eta_{\delta}-\D\ (v\cdot \nabla \eta_{\delta})$. Since $V(x,t)$ is the kernel of $e^{t\Delta}\mathbb{P}$, we see that 

\begin{align*}
v_{\delta}(x,t)
&=e^{t\Delta}v_{\delta}(0)+\int_{0}^{t}e^{(t-s)\Delta}\mathbb{P}f_{\delta}(s)\dd s\\
&=\int_{\mathbb{R}^{n}}\Gamma(x-y,t)v_{\delta}(y,0)\dd y+\int_{0}^{t}\int_{\mathbb{R}^{n}}V(x-y,t-s)f_{\delta}(y,s)\dd y\dd s.  
\end{align*}\\
The function $v_{\delta}(x,t)$ converges to $v(x,t)$ for $x\in \Omega$ and $t>0$ as $\delta\to0$. It is not difficult to see the convergence of the right-hand side. Since 

\begin{align*}
\int_{\mathbb{R}^{n}}\varphi(y)\cdot \nabla \eta_{\delta}(y)\dd y
&=-\int_{\{\delta<d(y)<2\delta\}}\D\ \varphi(y) \eta_{\delta}(y)\dd y
+\int_{\{d(y)=2\delta\}} \varphi(y)\cdot n(y)\dd y\\
&\to -\int_{\partial\Omega}\varphi(y)\cdot n(y)\dd y\quad \textrm{as}\ \delta\to0,\quad \varphi\in C^{1}_{0}(\mathbb{R}^{n}),
\end{align*}\\
we see that

\begin{align*}
\int_{\mathbb{R}^{n}}V(x-y,t-s)f_{\delta}(y,s)\dd y
&\to -\int_{\partial\Omega}V(x-y,t-s)(-\nabla v(y,s)+q(y,s)I)n(y)\dd H(y)\\
&=\int_{\partial\Omega}V(x-y,t-s)T(y,s)n(y)\dd H(y)
\end{align*}\\
for each $t>s>0$ by $\sum_{j}\partial_j V_{ij} =0$. Thus (3.3) is obtained by sending $\delta\to0$.

\noindent
(ii) The above observation implies that the formula (3.3) is extendable for complex time $t\in \Sigma_{\theta}$. In fact, by applying the Duhamel's principle on the segment $\gamma_{t}=\{s\in \gamma,\ |s|\leq |t|\ \}$ in the half line $\gamma=\{t=\arg{\theta}\}$ and sending $\delta\to0$, we have 

\begin{align*}
v(x,t)=\int_{\Omega}\Gamma(x-y,t)v_0(y)\dd y+\int_{\gamma_t}\int_{\partial\Omega}V(x-y,t-s)T(y,s)n(y)\dd H(y)\dd s   \tag{3.11}
\end{align*}\\
for $x\in \Omega$ and $t\in \gamma$. The Oseen tensor also satisfies the pointwise estimate

\begin{align*}
|V(x,t)|\leq \frac{C}{(|x|+|t|^{1/2})^{n}}\quad x\in \mathbb{R}^{n}, \ t\in \overline{\Sigma_{\theta}},
\end{align*}\\
with some constant $C=C(\theta)$ for $\theta\in [0,\pi/2)$. 
\end{rems}

\vspace{15pt}

%section4
\section{A large time estimate}

\vspace{10pt}

We prove (1.6) for $n\geq 3$ and $t>0$. Since the Stokes flow $v=S(t)v_0$ belongs to $D(A)$ for $t>0$, we are able to reduce (1.6) to the estimate by the graph-norm (4.4). We prove (4.4) by the Liouville Theorem (Theorem 1.1).

\vspace{10pt}

%prop4.1
\begin{lem}
There exists a constant $C$ such that 

\begin{align*}
\sup_{t>0}||S(t)v_0||_{L^{\infty}}\leq C||v_0||_{L^{\infty}}   \tag{4.1}
\end{align*}\\
holds for $v_0\in L^{\infty}_{\sigma}$ (i.e., the Stokes semigroup is a bounded semigroup on $L^{\infty}_{\sigma}$.)  
\end{lem}

\vspace{10pt}
We prove (4.1) for initial data with finite energy. 
\vspace{10pt}

%prop4.2
\begin{prop}
For $v_0\in L^{\infty}_{\sigma}\cap L^{2}$, the Stokes semigroup $S(t)v_0$ is bounded on $L^{\infty}_{\sigma}$ for all $t>0$, i.e., 

\begin{align*}
\sup_{t>0}||S(t)v_0||_{L^{\infty}}<\infty.   \tag{4.2}
\end{align*}
\end{prop}

\vspace{10pt}

\begin{proof}
Since $S(t)v_0$ is bounded on $L^{\infty}$ for $t\in (0 ,1]$ by (3.2), we consider the case $t\geq 1$. Let $A$ denote the Stokes operator on $L^{p}_{\sigma}$ with the domain $D(A)=W^{2,p}\cap W^{1,p}_{0}\cap L^{p}_{\sigma}$. Since $D(A)\subset W^{2,p}$ with continuous injection and the Stokes semigroup is bounded on $L^{p}_{\sigma}$, by the Sobolev embedding for $p>n$, it follows that 

\begin{align*}
||S(t)v_0||_{L^{\infty}}
\leq C||S(t)v_0||_{W^{1,p}}
\leq C'(||S(t)v_0||_{L^{p}}+||AS(t)v_0||_{L^{p}})
\leq C''||v_0||_{L^{p}},\quad t\geq 1.
\end{align*}\\
Thus (4.2) holds.
\end{proof}

\vspace{10pt}

%lem4.3
\begin{prop}
The estimate (4.1) holds for all $v_0\in L^{\infty}_{\sigma}\cap L^{2}$.
\end{prop}

\vspace{5pt}

\begin{proof}[Proof of Lemma 4.1]
For $v_0\in L^{\infty}_{\sigma}$, we take a sequence $\{v_{0,m}\}\subset C^{\infty}_{c,\sigma}$ such that 

\begin{equation*}
\begin{aligned}
||v_{0,m}||_{L^{\infty}}\leq C||v_{0}||_{L^{\infty}}, \\
v_{0,m}\to v_{0}\quad \textrm{a.e. in}\ \Omega,
\end{aligned}
\tag{4.3}
\end{equation*}\\
by \cite[Lemma 5.1]{AG2}. Since $S(t)v_{0,m}$ converges to $S(t)v_0$ locally uniformly in $\overline{\Omega}\times (0,\infty)$ \cite{AG2} and (4.1) holds for $v_{0,m}$ by Proposition 4.3, (4.1) is extendable for all $v_0\in L^{\infty}_{\sigma}$. 
\end{proof}

\vspace{10pt}
Proposition 4.3 follows from:

\vspace{10pt}

%prop4.4
\begin{prop}
There exists a constant $C$ such that
 
\begin{align*}
\sup_{t>0}||S(t)v_0||_{L^{\infty}}(t)\leq C||v_0||_{D(A)}  \tag{4.4}
\end{align*}\\
holds for $v_0\in L^{\infty}_{\sigma}\cap L^{2}$ satisfying $v_0\in D(A)$, where $D(A)$ is the domain of the Stokes operator on $L^{\infty}_{\sigma}$.
\end{prop}

\vspace{10pt}

\begin{proof}[Proof of Proposition 4.3]
Since $S(t)$ is bounded on $L^{\infty}_{\sigma}$ for $t\in (0,1]$, we consider the case $t\geq 1$. Since $S(1)v_0\in L^{\infty}_{\sigma}\cap L^{2}\cap D(A)$ for $v_0\in L^{\infty}_{\sigma}\cap L^{2}$ and 

\begin{align*}
||S(1)v_0||_{D(A)}\leq C||v_0||_{L^{\infty}}
\end{align*}\\
by (3.2), it follows from (4.4) that 

\begin{align*}
||S(t)v_0||_{L^{\infty}}=||S(t-1)S(1)v_0||_{L^{\infty}}\leq C||S(1)v_0||_{D(A)}
\leq C'||v_0||_{L^{\infty}}.
\end{align*}\\
Thus (4.1) holds for all $v_0\in L^{\infty}_{\sigma}\cap L^{2}$.
\end{proof}

\vspace{10pt}

We now prove (4.4) by the Liouville theorem.

\vspace{10pt}

\begin{proof}[Proof of Proposition 4.4]
We argue by a contradiction. Suppose on the contrary that (4.4) were false. Then, for any $m\geq 1$ there exists $\tilde{v}_{0,m}\in L^{\infty}_{\sigma}\cap L^{2}\cap D(A)$ such that 

\begin{align*}
M_m:=\sup_{t>0}||\tilde{v}_{m}||_{L^{\infty}}(t)> m||\tilde{v}_{0,m}||_{D(A)}
\end{align*}\\
for $\tilde{v}_m=S(t)\tilde{v}_{0,m}$. We set $v_{m}=\tilde{v}_{m}/M_m$ so that 

\begin{align*}
\sup_{t>0}||v_m||_{L^{\infty}}(t)=1,\quad ||v_{0,m}||_{D(A)}< \frac{1}{m}.
\end{align*}\\
We take the associated pressure $q_m$ satisfying 

\begin{align*}
\int_{\Omega_0}q_m(x,t)\dd x=0,
\end{align*}\\
for $\Omega_0=B_0(R_0)\cap \Omega$ and $R_0>0$ such that $\Omega^{c}\subset B_0(R_0)$, and denote the stress tensor by $T_m=\nabla v_m+\nabla^{T} v_m-q_mI$. Then, by Lemma 3.3, there exists a constant $C$, independent of $m$ such that 

\begin{align*}
\sup_{t>0}||T_m||_{L^{\infty}(\partial\Omega)}(t)
\leq C.
\end{align*}\\
We take $t_m\in (0,\infty)$ such that $||v_m||_{\infty}(t_m)\geq  1/2$. By (3.2), we may assume that $t_m\to\infty$. We take $x_m\in \Omega$ such that $|v_m(x_m,t_m)|\geq 1/4$.\\

\noindent
\textit{Case 1.} $\overline{\lim}_{m\to\infty}d(x_m)=\infty$. We may assume that $\lim_{m\to\infty}d(x_m)=\infty$ by choosing a subsequence. By Lemmas 3.2 and 3.3, it follows that 

\begin{align*}
\frac{1}{4}
\leq |v_m(x_m,t_m)|
&\leq ||v_{0,m}||_{L^{\infty}}+\frac{C}{d(x_m)^{n-2}}\sup_{0<s\leq t_m}||T_m||_{L^{\infty}(\partial\Omega)}(s)\\
&\leq \frac{1}{m}+\frac{C'}{d(x_m)^{n-2}}\to 0\quad \textrm{as}\ m\to\infty.
\end{align*}\\
Thus Case 1 does not occur.\\

\noindent
\textit{Case 2.} $\overline{\lim}_{m\to\infty}d(x_m)<\infty$. We may assume that $x_m\to x_{\infty}\in \overline{\Omega}$ by choosing a subsequence. We set 
\begin{align*}
u_m(x,t)=v_m(x,t+t_m),\quad p_m(x,t)=q_m(x,t+t_m),
\end{align*}\\
so that $(u_m,p_m)$ is a solution of the Stokes equations in $\Omega\times (-t_m,0]$. By Lemmas 3.2 and 3.3, it follows that 
\begin{align*}
|u_m(x,t)|\leq \frac{1}{m}+\frac{C}{d(x)^{n-2}}\quad x\in \Omega,\ t\in (-t_m,0]. 
\end{align*}\\
Since $u_m$ is bounded in $\Omega\times (-t_m,0]$, it follows from Proposition 3.1 that $\partial_t^{s}\partial_x^{k}u_m$ are bounded on $\overline{\Omega}\times (-T,0]$ for $2s+|k|\leq 2$ and each fixed $T>0$. Thus, there exists a subsequence such that $u_m$ converges to a limit $u$ locally uniformly in $\overline{\Omega}\times (-\infty,0]$. By sending $m\to\infty$, the limit satisfies

\begin{align*}
|u(x,t)|\leq \frac{C}{d(x)^{n-2}}\quad x\in \Omega,\ t\in (-\infty,0]. \tag{4.5}
\end{align*}\\
We see that the limit $u$ is an ancient solution to (1.1). We take $\varphi\in C^{2,1}_{c}(\overline{\Omega}\times (-\infty,0])$ satisfying $\D\ \varphi=0$ in $\Omega\times (-\infty,0)$ and $\varphi=0$ on $\partial\Omega\times (-\infty,0)\cup \Omega\times \{t=0\}$. Since $u_m$ satisfies (1.1) in $\Omega\times (-t_m,0]$ and $\varphi$ is supported in $\overline{\Omega}\times (-t_m,0]$ for sufficiently large $m$, by multiplying $\varphi$ by (1.1) and integration by parts, it follows that

\begin{align*}
\int_{-t_m}^{0}\int_{\Omega}u_m\cdot (\partial_t \varphi+\Delta \varphi)\dd x\dd t=0.
\end{align*}\\
Sending $m\to\infty$ implies that $u$ is an ancient solution to (1.1). Since $u\in L^{\infty}(-\infty,0; L^{p})$ for $p\in (n/(n-2),\infty)$ by (4.5), applying the Liouville theorem (Theorem 1.1) implies that $u\equiv 0$. This contradicts $|u(x_\infty,0)|\geq 1/4$. Thus Case 2 does not occur.

We reached a contradiction. The proof is now complete.
\end{proof}

\vspace{10pt}

\section{Extensions to complex time}

\vspace{10pt}
We prove Theorem 1.2. We first reduce the estimate (1.6) in the sector $\Sigma_{\theta}$ to that on the half lines $\{\arg{t}=\pm\theta\}$ by the maximum principle. In the subsequent section, we prove (1.6) on the half lines by applying the Liouville theorem in a sector (Theorem 2.7). After the proof of Theorem 1.2, we prove Theorem1.3.

\vspace{10pt}

\subsection{A maximum in a sector}

%lem5.1
\begin{lem}
For $\theta\in [0,\pi/2)$, there exists a constant $C$ such that 

\begin{align*}
\sup_{t\in \Sigma_{\theta}}||S(t)v_0||_{L^{\infty}}\leq C||v_0||_{L^{\infty}}  \tag{5.1}
\end{align*}\\
holds for $v_0\in L^{\infty}_{\sigma}$.
\end{lem}

\vspace{10pt}
\begin{proof}[Proof of Theorem 1.2]
The assertion follows from Lemma 5.1.
\end{proof}
\vspace{10pt}

\vspace{10pt}
The goal of this subsection is to prove Lemma 5.1 by using: 
\vspace{10pt}

%lem5.2
\begin{lem}
There exists a constant $C$ such that 

\begin{align*}
\sup\left\{ ||S(t)v_0||_{L^{\infty}}\ |\ t\in \mathbb{C}\backslash \{0\},\ \arg{t}=\pm\theta \right\}\leq C||v_0||_{L^{\infty}}  \tag{5.2}
\end{align*}\\
for $v_0\in L^{\infty}_{\sigma}\cap L^{2}$.
\end{lem}

\vspace{10pt}

%prop5.3
\begin{prop}
Let $v_0\in L^{\infty}_{\sigma}\cap L^{2}$. Suppose that there exists a constant $M>0$ such that 

\begin{align*}
\sup\left\{||S(t)v_0||_{L^{\infty}}\ \middle|\ t\in \mathbb{C}\backslash \{0\},\ \arg{t}=\theta   \right\}\leq M.   \tag{5.3}
\end{align*}\\
Then, 

\begin{align*}
\sup\left\{||S(t)v_0||_{L^{\infty}}\ \middle|\ t\in \Sigma_{\theta},\ 0\leq \arg{t}\leq \theta   \right\}\leq \max\{M, C||v_0||_{L^{\infty}}\}.  \tag{5.4}
\end{align*}\\
with some constant $C$, independent of $v_0$. 
\end{prop}

\vspace{10pt}

%lem5.2
\begin{prop}
The estimate (5.1) holds for $v_0\in L^{\infty}_{\sigma}\cap L^{2}$.
\end{prop}

\vspace{5pt}

\begin{proof}
Since (5.3) holds for $M=C||v_0||_{L^{\infty}}$ by (5.2), we apply Proposition 5.3 and estimate

\begin{align*}
\sup\{||S(t)v_0||_{L^{\infty}}|\ t\in \mathbb{C}\backslash \{0\},\ 0\leq \arg{t}\leq \theta  \}\leq C||v_0||_{L^{\infty}}.
\end{align*}\\
Since the assertion of Proposition 5.3 is also valid for the half line $\{\arg{t}=-\theta\}$, we are able to estimate the sup-norm in $\{-\theta\leq \arg{t}\leq 0\}$ and obtain (5.1).
\end{proof}

\vspace{5pt}
\begin{proof}[Proof of Lemma 5.1]
The estimate (5.1) is extendable for $v_0\in L^{\infty}_{\sigma}$ by the approximation (4.3) and Proposition 5.4.
\end{proof}

\vspace{10pt}

We prove Proposition 5.3. We use a decay property in order to estimate a maximum in a sector $\Sigma_{\theta}$.

\vspace{10pt}

%prop5.4
\begin{prop}
For $v_0\in L^{\infty}_{\sigma}\cap L^{2}$, we have 

\begin{align*}
\lim_{T\to\infty}\sup\left\{||S(t)v_0||_{L^{\infty}}\ \middle|\
\textrm{Re}\ t=T,\ 0\leq \arg{t}\leq \theta \right\}=0.  \tag{5.5}
\end{align*}
\end{prop}

\vspace{5pt}

\begin{proof}
We show that 

\begin{align*}
\lim_{T\to\infty}\sup\left\{||S(t)v_0||_{L^{p}}\ \middle|\
\textrm{Re}\ t=T,\ 0\leq \arg{t}\leq \theta' \right\}=0,   \tag{5.6}
\end{align*}\\
for $\theta'\in [0,\pi/2)$. Let $A$ denote the Stokes operator on $L^{p}_{\sigma}$ for $p\in (1,\infty)$ with the domain $D(A)=W^{2,p}\cap W^{1,p}_{0}\cap L^{p}_{\sigma}$. Since the range of the Stokes operator $A$ is dense in $L^{p}_{\sigma}$ \cite{GS}, \cite{BV}, we have $\lim_{t\to\infty}||S(t)v_0||_{L^{p}}=0$. We take $\theta''\in (\theta',\pi/2)$. For $t\in \mathbb{C}$ satisfying $\textrm{Re}\ t=T$ and $0\leq \arg{t}\leq \theta'$, there exists $t_1\in [\kappa T,T]$ such that $t=t_1+|t-t_1|e^{i\theta''}$ for $\kappa=(1-\tan{\theta'}/\tan{\theta''})$. Since $S(t)v_0$ is bounded on the half line $\{\arg{t}=\theta''\}$, it follows that 

\begin{align*}
||S(t)v_0||_{L^{p}}=||S(|t-t_1|e^{i\theta''})S(t_1)v_0||_{L^{p}}
\leq C\sup_{\kappa T\leq t_1\leq T}||S(t_1)v_0||_{L^{p}}.
\end{align*}\\
By taking a supremum for $t\in \{\textrm{Re}t=T,\ 0\leq \arg{t}\leq \theta'\}$ and sending $T\to\infty$, we see that (5.6) holds. Since 

\begin{align*}
||S(t)v_0||_{L^{\infty}}
\leq C||S(t)v_0||_{W^{1,p}}
&\leq C'(||S(t)v_0||_{L^{p}}+||A S(t)v_0||_{L^{p}}  ) \\
&\leq C''(||S(t)v_0||_{L^{p}}+||S(t-1)v_0||_{L^{p}} ),
\end{align*}\\
for $p>n$ as in the proof of Proposition 4.2, (5.5) follows from (5.6).
\end{proof}

\vspace{5pt}
\begin{proof}[Proof of Proposition 5.3]
We set $A_{\varepsilon,T}=\{t\in \mathbb{C}\backslash \{0\}\ |\ \varepsilon< \textrm{Re}\ t\leq  T,\ 0\leq \arg{t}\leq \theta  \}$. Since $S(t)v_0$ is holomorphic and continuous in $\overline{A_{\varepsilon,T}}$, we estimate

\begin{align*}
\sup_{t\in A_{\varepsilon,T}}||S(t)v_0||_{L^{\infty}}\leq \sup_{t\in \partial A_{\varepsilon,T}}||S(t)v_0||_{L^{\infty}}
\end{align*}\\
by the maximum principle. By Proposition 3.1 (ii) and Lemma 4.1, we estimate 

\begin{align*}
||S(t)v_0||_{L^{\infty}}\leq C||v_0||_{L^{\infty}}
\end{align*}\\
for $t\in \{0<\textrm{Re}\ t\leq \varepsilon,\ 0\leq \arg{t}\leq \theta \}\cup \{t>0\}$. It follows that  

\begin{align*}
&\sup_{t\in A_{0,T}}||S(t)v_0||_{L^{\infty}}
\leq  \max\{M, C||v_0||_{L^{\infty}}\}+R_{T},\\
&R_T=\sup\left\{||S(t)v_0||_{L^{\infty}}\ \middle|\
\textrm{Re}\ t=T,\ 0\leq \arg{t}\leq \theta \right\}.
\end{align*}\\
Since $R_T\to0$ as $T\to\infty$ by (5.5), we obtain (5.4).
\end{proof}

\vspace{10pt}

\vspace{10pt}

\subsection{Boundedness on half lines}

\vspace{10pt}
It remains to show Lemma 5.2. 

\vspace{10pt}

%prop5.7
\begin{prop}
There exists a constant $C$ such that 

\begin{align*}
\sup\left\{||S(t)v_0||_{L^{\infty}}\ |\ t\in \mathbb{C}\backslash\{0\},\ \arg{t}=\pm\theta \right\}\leq C||v_0||_{D(A)}  \tag{5.8}
\end{align*}\\
holds for $v_0\in L^{\infty}_{\sigma}\cap L^{2}\cap D(A)$. 
\end{prop}

\vspace{10pt}

\begin{proof}[Proof of Lemma 5.2]
We take $t\in \mathbb{C}\backslash \{0\}$ such that $\arg{t}=\theta$. Since (5.2) holds for $|t|\leq 1$ by Proposition 3.1 (ii), we consider $|t|\geq 1$. By the semigroup property $S(t)v_0=S((|t|-1)e^{i\theta})S(e^{i\theta})v_0$ and $||AS(e^{i\theta})v_0||_{L^{\infty}}\leq C||v_0||_{L^{\infty}}$, it follows from (5.8) that 

\begin{align*}
||S(t)v_0||_{L^{\infty}}
=||S((|t|-1)e^{i\theta})S(e^{i\theta})v_0||_{L^{\infty}}
\leq C||S(e^{i\theta})v_0||_{D(A)}
\leq C'||v_0||_{L^{\infty}}.
\end{align*}\\
Thus (5.2) holds on $\{\arg{t}=\theta\}$. By the same way, we obtain (5.2) on $\{\arg{t}=-\theta\}$.
\end{proof}

\vspace{10pt}

The estimates used to prove (1.6) for $t>0$ are extendable for complex time.

\vspace{10pt}

%prop5.7
\begin{prop}
(i) Set $\gamma=\{t\in \mathbb{C}\backslash \{0\}\ |\ \arg{t}=\theta \}$ for $\theta\in [0,\pi/2)$. For $T>0$, there exists a constant $C$ such that 

\begin{align*}
\sup\left\{|t|^{\frac{|k|}{2}+s} ||\partial_t^{s}\partial_x^{k}S(t)v_0||_{L^{\infty}}\middle|\ t\in \gamma,\ \textrm{Re}\ t\leq T \right\}\leq C||v_0||_{L^{\infty}}  \tag{5.9}
\end{align*}\\
holds for $v_0\in L^{\infty}_{\sigma}$ and $2s+|k|\leq 2$.

\noindent 
(ii) For $n\geq 3$, there exists a constant $C$ such that 

\begin{align*}
|v(x,t)|\leq C\left(||v_0||_{L^{\infty}}+\frac{1}{d(x)^{n-2}}\sup_{s\in \gamma_{t}}||T||_{L^{\infty}(\partial\Omega)}(s) \right)   \tag{5.10}
\end{align*}\\
for $x\in \Omega$ and $t\in \gamma$.

\noindent
(iii) There exists a constant $C$ such that 

\begin{align*}
\sup_{s\in \gamma_t}||T||_{L^{\infty}(\partial\Omega)}(s)
\leq C\left(||v_0||_{D(A)}+\sup_{s\in \gamma_t}||v||_{L^{\infty}}(s) \right)  \tag{5.11}
\end{align*}\\
holds for all $t\in \gamma$ and $v_0\in D(A)$, with associated pressure satisfying (3.5).
\end{prop}

\vspace{10pt}

\begin{proof}
The estimate (5.9) follows from Proposition 3.1. Since the formula (3.11) holds for complex time $t\in \Sigma_{\theta}$ as in Remarks 3.5 (ii), the pointwise estimate (5.10) holds. The estimate (5.11) follows from (5.9) and (3.9). 
\end{proof}

\vspace{10pt}

We now complete:
\vspace{10pt}

\begin{proof}[Proof of Proposition 5.6]
We prove (5.8) on the half line $\gamma=\{\arg{t}=\theta\}$. By the same way, we are able to prove (5.8) on $\{\arg{t}=-\theta\}$. We argue by a contradiction. Suppose that (5.8) were false. Then, for $m\geq 1$ there exisits $v_{0,m}\in L^{\infty}_{\sigma}\cap L^{2}\cap D(A)$ such that 

\begin{align*}
\sup_{t\in \gamma}||v_{m}||_{L^{\infty}}(t)=1,\quad ||v_{0,m}||_{D(A)}<\frac{1}{m},
\end{align*}\\
for $v_m=S(t)v_{0,m}$. Since $||v_m||_{L^{\infty}}$ is bounded on $\{t>0\}$ by (4.1), we apply Proposition 5.3 to estimate 

\begin{align*}
\sup\left\{||v_m||_{L^{\infty}}(t)\ \middle|\ t\in \mathbb{C}\backslash\{0\},\ 0\leq  \arg{t}\leq  \theta\right\}\leq C,  \tag{5.12}
\end{align*}\\
with some constant $C$, independent of $m$. We take the associated pressure $q_m$ satisfying (3.5) and set $T_m=\nabla v_m+\nabla^{T} v_m-q_mI$. By (5.11) and (5.12), $||T_m||_{L^{\infty}(\partial\Omega)}(t)$ is uniformly bounded on $\{\arg{t}=\theta'\}$ for all $\theta'\in [0,\theta]$, i.e., 

\begin{align*}
\sup\left\{||T_m||_{L^{\infty}(\partial\Omega)}(t)\ \middle|\ t\in \mathbb{C}\backslash\{0\},\ 0\leq \arg{t}\leq \theta\right\}
\leq C'.  \tag{5.13}
\end{align*}\\
We take $t_m\in \gamma$ such that $||v_m||_{L^{\infty}}(t_m)\geq 1/2$. We may assume that $|t_m|\to\infty$ by (5.9). We take $x_m\in \Omega$ such that $|v_m(x_m,t_m)|\geq 1/4$.\\

\noindent 
\textit{Case 1.}\ $\overline{\lim}_{m\to\infty}d(x_m)=\infty$. We may assume that ${\lim}_{m\to\infty}d(x_m)=\infty$. By (5.10), it follows that

\begin{align*}
\frac{1}{4}\leq |v(x_m,t_m)|\leq C\left(\frac{1}{m}+\frac{1}{d(x_m)^{n-2}}\right)\to0
\quad \textrm{as}\ m\to\infty.
\end{align*}\\
Thus Case 1 does not occur.\\

\noindent 
\textit{Case 2.}\ $\overline{\lim}_{m\to\infty}d(x_m)<\infty$. We may assume that $x_m\to x_{\infty}\in \overline{\Omega}$ by choosing a subsequence. We set 

\begin{align*}
u_m(x,t)=v_m(x,t+t_m),\quad p_m(x,t)=q_m(x,t+t_m)
\end{align*}\\
so that $(u_m,p_m)$ satisfies the Stokes equations (1.1) in $\Omega\times \Lambda_m$ for 

\begin{align*}
\Lambda_m=\left\{t\in \mathbb{C}\backslash \{0\}\ \middle|\ t\neq -t_m,\ \pi-\theta\leq \arg{t}\leq -\frac{\pi}{2},\ \textrm{Im}\ t\geq -|t_m|\sin\theta     \right\}.
\end{align*}\\
The domain $\Lambda_m$ approaches to the sector 

\begin{align*}
\Lambda=\left\{t\in \mathbb{C}\backslash\{0\}\ \middle|\ \pi-\theta\leq \arg{t}\leq -\frac{\pi}{2}\right\}.
\end{align*}\\
It follows from (5.12) and (5.13) that

\begin{align*}
\sup\left\{||u_m||_{L^{\infty}(\Omega)}(t)+||T_m||_{L^{\infty}(\partial\Omega)}(t)\ \middle|\ t\in \Lambda_m\right\}\leq C. \tag{5.14} 
\end{align*}\\
We apply (5.10) to estimate

\begin{align*}
|u_m(x,t)|\leq C\left(\frac{1}{m}+\frac{1}{d(x)^{n-2}}\right)  \tag{5.15}
\end{align*}\\
for $x\in \Omega$ and $t\in \Lambda_m$. Since $\partial_t^{s}\partial_x^{k}u_m$ for $2s+|k|\leq 2$ are uniformly bounded for each bounded domain in $\overline{\Omega}\times \overline{\Lambda}$ by (5.9), there exists a subsequence such that $u_m$ converges to a limit $u$ locally uniformly in $\overline{\Omega}\times \overline{\Lambda}$. For each $T>0$, we set 
$I_T=\{t\in \Lambda\ |\ -T\leq \textrm{Re}\ t\leq 0,\ \textrm{Im}\ t=-T\tan{\theta}\}$. Since $I_T\subset \Lambda_m$ for sufficiently large $m$, $(u_m,p_m)$ satisfies (1.1) in $\Omega\times I_T$. We take an arbitrary $\varphi\in C^{2,1}_{c}(\overline{\Omega}\times [-T,0])$ satisfying $\D\ \varphi=0$ in $\Omega\times [-T,0]$ and $\varphi=0$ on $\partial\Omega\times (-T,0)\cup \Omega\times \{t=0\}$. By multiplying $\varphi$ by (1.1) in $\Omega\times I_T$ and integration by parts, it follows that 

\begin{align*}
\int_{-T}^{0}\int_{\Omega}u_m(x,\alpha+i\beta)(\partial_\alpha\varphi+\Delta \varphi)\dd x\dd \alpha=-\int_{\Omega}u_m(x,-T+i\beta)\varphi(x,-T+i\beta)\dd x,\quad \beta=-T\tan{\theta}.
\end{align*}\\
Sending $m\to\infty$ implies that 

\begin{align*}
\int_{-T}^{0}\int_{\Omega}u(x,\alpha+i\beta)(\partial_\alpha\varphi+\Delta \varphi)\dd x\dd \alpha=-\int_{\Omega}u(x,-T+i\beta)\varphi(x,-T+i\beta)\dd x.
\end{align*}\\
Thus, the limit $u$ is an ancient solution in $\Omega\times \Lambda$. Since $u\in L^{\infty}(\Lambda; L^{p})$ for $p\in (n/(n-2),\infty)$ by (5.15), applying the Liouville theorem (Theorem 2.7) implies $u\equiv 0$. This contradicts $|u(x_{\infty},0)|\geq 1/4$. Thus Case 2 does not occur.\\
We reached a contradiction. The proof is now complete.
\end{proof}

\vspace{10pt}

\begin{proof}[Proof of Theorem 1.3]
For $n=2$ and solutions $v=S(t)v_0$ with zero net force, the estimate

\begin{align*}
|v(x,t)|\leq ||v_0||_{L^{\infty}(\Omega)}+\frac{C}{|x|}\sup_{0<s\leq t}||T||_{L^{\infty}(\partial\Omega)}(s)  \tag{5.16}
\end{align*}\\
holds for $|x|\geq 2R_0$ and $t>0$ with $R_0\geq \textrm{diam}\ \Omega^{c}$. By using (5.16) and the Liouville theorem (Theorem 1.1), we are able to prove (4.4) by the same way as $n\geq 3$ and (4.1) holds. The complex time case is parallel. Thus Lemmas 4.1 and 5.1 hold for $n=2$ and $v=S(t)v_0$ with zero net force.
\end{proof}

\vspace{10pt}
\begin{rems}
(i) Theorems 1.2 and 1.3 imply the large time estimate

\begin{align*}
\sup_{t>0}\left\{||S(t)v_0||_{L^{\infty}}+t||AS(t)v_0||_{L^{\infty}}\right\}\leq C||v_0||_{L^{\infty}}  \tag{5.17}
\end{align*}\\
for $v_0\in L^{\infty}_{\sigma}$. The estimate (5.17) follows from a resolvent estimate. Since the resolvent is represented by $S(t)$ by the Laplace transform and the integral path is replaced from $(0,\infty)$ to the half line $\gamma=\{\arg{t}=\theta\}$ for $\theta\in (0,\pi/2)$, we see that

\begin{align*}
(\lambda-A)^{-1}f
=\int_{0}^{\infty}e^{-\lambda t}S(t)f\dd t =\int_{\gamma}e^{-\lambda t}S(t)f\dd t.
\end{align*}\\
Since $S(t)$ is bounded in the sector $\Sigma_{\theta}$, we obtain 

\begin{align*}
||(\lambda-A)^{-1}f||_{L^{\infty}}\leq \frac{C}{|\lambda|}||f||_{L^{\infty}}\quad \lambda\in \Sigma_{\theta'+\pi/2}  \tag{5.18}
\end{align*}\\
with some constant $C=C_{\theta'}$ for $\theta'\in [0,\pi/2)$. The estimate (5.18) implies (5.17).

\noindent
(ii) The spatial derivative estimate  

\begin{align*}
||\nabla S(t)v_0||_{L^{\infty}}\leq \frac{C}{t^{1/2}}||v_0||_{L^{\infty}}\quad 0<t\leq T   \tag{5.19}
\end{align*}\\
holds for each $T>0$ \cite{AG2}. We are not able to take $T=\infty$ in (5.19). To see this, we recall the decay estimate 

\begin{align*}
||\nabla S(t)v_0||_{L^{p}}\leq \frac{C}{t^{n/2(1/q-1/p)+1/2 }}||v_0||_{L^{q}}\quad t>0, \tag{5.20}
\end{align*}\\
for $v_0\in L^{q}_{\sigma}$ and $1<q\leq p\leq n$. See \cite{Iwashita89} for $n\geq 3$, \cite{DS1} for $n=2$ and \cite{MS97}. It is known that the condition $p\leq n$ is optimal \cite{MS97}, \cite[Corollary 2.4]{H11} in the sense that (5.20) for $p>n$ is not valid for all $t\geq 1$ and $v_0\in L^{q}_{\sigma}$. If (5.19) were true for all $t>0$, by the semigroup property and the decay estimate 

\begin{align*}
||S(t)v_0||_{L^{\infty}}\leq \frac{C}{t^{n/(2q)}}||v_0||_{L^{q}}\quad t>0,
\end{align*}\\
proved in \cite{Iwashita89},\cite{Chen93} for $n\geq 3$, \cite{DS2} for $n=2$, we would obtain (5.20) for $p=\infty$.

\noindent 
(iii) Theorem 1.3 improves the pointwise estimates of the two-dimensional Navier-Stokes flows for rotationally symmetric initial data around a unit disk $\Omega^{c}$. Let $u$ be a global-in-time solution of the Navier-Stokes equations for initial data $u_0\in L^{2}_{\sigma}\cap L^{1}\cap W^{2-2/q,q}_{0}(\Omega)$ for $q\in (1,4/3]$. It is proved in \cite[Theorem 5.8]{HeMiyakawa06} that if $u_0$ is $D_{m+2}$-covariant for some $m\geq 0$ and satisfies $u_0\in W^{1,2}_{0}(\Omega)$, $(1+|x|)^{m+3}|u_0(x)|\in L^{\infty}(\Omega)$ and $(1+|x|)^{m+1}|u_0(x)|\in L^{1}(\Omega)$, then $u$ is $D_{m+2}$-covariant and satisfies the pointwise estimates

\begin{align*}
|u(x,t)|\lesssim 
\begin{cases}
& |x|^{-(m+3)}\quad |x|\geq 2,\ t>0,\\
& t^{-(m+3)/2}\quad x\in \Omega,\ t>0.
\end{cases}
\tag{5.21}
\end{align*}\\
The estimate (5.21) is obtained from the representation formula of the Navier-Stokes flows. Although the right-hand side is unbounded at $x=0$ and $t=0$, respectively, by estimating the integral form of $u$ by using (1.6), we are able to show that $u$ is bounded in $\Omega\times (0,\infty)$. Hence (5.21) is improved to 

\begin{align*}
|u(x,t)|\lesssim 
\begin{cases}
& (1+|x|)^{-(m+3)},\\
& (1+t)^{-(m+3)/2}\quad x\in \Omega,\ t>0,
\end{cases}
\end{align*}\\
as noted in \cite[p.1546, Remarks (ii)]{HeMiyakawa06}.

\end{rems}

\vspace{10pt}
\section*{Acknowledgements}
The author is grateful to Professor Toshiaki Hishida for informing him of the paper \cite{H11} on the optimal exponent of the decay estimate (5.20). This work is partially supported by JSPS through the Grant-in-aid for Young Scientist (B) 17K14217, Scientific Research (B) 17H02853 and Osaka City University Strategic Research Grant 2018 for young researchers.

\vspace{10pt}

%ref
\bibliographystyle{plain}
\bibliography{ref}

\end{document}